\documentclass[11pt,twoside]{article}
\usepackage{euscript} 
\usepackage{mathrsfs}
\usepackage{bm} 
\usepackage{times}
\usepackage[tmargin=25mm,bmargin=25mm,lmargin=25mm,rmargin=25mm]{geometry}
\fontsize{12}{14}\selectfont

\date{} 

\begin{document} 


\centerline{}

\centerline {\Large{\bf A Study of the I - functions of two variables}} 

\centerline{} 

\centerline{\bf{Shantha Kumari.K.}}
\centerline{\footnotesize Research Scholar, SCSVMV, Kanchipuram,}
\centerline{\footnotesize and  Department of Mathematics,}
\centerline{\footnotesize P.A.College of Engineering, Mangalore,}
 \centerline{\footnotesize Karnataka, INDIA, E-mail: skk\_abh@rediffmail.com}

\centerline{}

\centerline{\bf{Vasudevan Nambisan T.M.}}
\centerline{\footnotesize Department of Mathematics, College of Engineering,}
\centerline{\footnotesize Trikaripur, Kerala, INDIA}
\centerline{\footnotesize Email:tmvasudevannambisan@yahoo.com}

\centerline{}

\centerline{\bf{Arjun K. Rathie}}
\centerline{\footnotesize Department of Mathematics,}
\centerline{\footnotesize School of Mathematical and Physical sciences,}
\centerline{\footnotesize Central University of Kerala, Kasaragod,INDIA.}
\centerline{\footnotesize Email: akrathie@gmail.com}

 \centerline{}
 \centerline{}

\begin{abstract}
 In our present investigation we propose to study and develop the I-function of two variables analogous to the I-function of one variable introduced  and studied by one of the  authors[24]. The conditions for convergence, series representation, behaviour for small values, elementary properties, transformation formulas and some special cases for the I-function of two variables are also discussed.
\end{abstract}
 


\normalsize
\section{Introduction}
\hspace{20 pt} The well known H-function of one variable defined by Fox [11] was in fact first \linebreak introduced by Pincherle in 1888 [10, section 1.19], but when Fox investigated and proved the H-function as a symmetric Fourier kernel to Meijers's G-function [10, p.206-222], the interest of most researchers , mathematicians and statisticians in this function has increased. By this fact, the H-function is often called Fox's H-function. Later on,  in 1964, Braaksma[6] studied and developed this function in a resonable good manner by finding  its asymptotic \linebreak behaviour etc. Moreover, the H-function is recognized as a generalization of both\linebreak  Meijer's G-function and  Wright's generalized hypergeometric functions[10,p.183]. \linebreak Numerous research papers by the researchers, mathematicians and statisticians related to \linebreak Fox's H-function, its properties and applications have been appeared in the literature. For more information about the results, we refer [17, 19] to the readers.

 \par  Recently  good deal of progress has been done in the direction of generalizing the Fox's H-function and Meijer's G-function. Frankly speaking, the well known H-function of one variable, introduced by Fox[11] and studied by Braaksma[6] contains as particular cases most of the  commonly used  special  functions of applied mathematics, but it does not \linebreak  contain some of the important functions such as the Riemann zeta functions, \linebreak polylogarithms etc. By demonstrating several examples of functions  which are not \linebreak included in the Fox's H-function, in 1997,  Rathie[24] introduced a new function in the\linebreak  literature namely the  ``I-function"  which is useful in Mathematics, physics and other branches of applied  mathematics. The newly defined function contains the polylogarithms, the exact partition of Gaussian free energy model from statistical mechanics, Feynmann integrals and functions useful in testing hypothesis from statistics as special cases.

 \par  Very recently, the I-function introduced by Rathie [24] has found useful applications in \textbf{wireless communications}. It is not out of place to mention here that, in two  of papers, Ansari, et. al. [2,3] have successfully developed an efficient Mathematica$^@$ implementation of the I-function, in order to give numerical results of their research. In one of their useful and very interesting papers, Ansari, et.al. [3] derived novel closed form expressions for the PDF and CDF of the sum of independently but not identically distributed (i.n.i.d.) gamma or equivalently squared Nakagami- m random variables in the case of both integer-order as well as non integer-order fading figure parameters. They have expressed their results in terms of Fox's H - function and Rathie's I-function.

\par The vast popularity and immense usefulness of generalized hypergeometric functions in one variable, inspired and stimulated a number of research workers to the study of hypergeometric functions involving two variables. Serious and significant study of the functions of two variables began with the introduction of $F_1$ , $F_2$ , $F_3$ and $F_4$ by Appell [4] and their confluent form by Humbert[16]. These  functions were further generalized by Kamp$\acute{e}$ de F$\acute{e}$riet by means of a function popularly known as the Kamp$\acute{e}$ de F$\acute{e}$riet function.

\par Appell and Kamp$\acute{e}$ de F$\acute{e}$riet [4]  studied the functions of two variables in details and recorded all the useful and important results concerning these functions in their famous work . Later on Srivastava and Daoust[29] studied the Kamp$\acute{e}$ de F$\acute{e}$riet function in a more general form. The hypergeometric functions of two variables have attracted attention of \linebreak eminent mathematicians. Thus Bailey,W.N.[5], Burchhal,J.L. and  Chaundy,T.W.[8,9], \linebreak Erdelyi,A.[10] etc. have contributed a lot of their developement and progress.

\par Further generalization of these functions of two variables were introduced almost \linebreak simultaneously by Sharma [27,28] and Agarwal [1]. These two functions, infact are the same except in their notational representations. Agarwal has given a method of \linebreak estimating this function for large values of the variable and has also obtained a formal pair of \linebreak unsymmetrical Fourier kernals for it. While Sharma[27,28]  has discussed various solutions of the partial differential equations satisfied by his function. He has also evaluated a number of interesting integrals and established several properties of the function.

\par Since then a number of mathematicians notably Pathak[22], Bora and Kalla [7],  Verma [31] and several others have studied functions of two variables which are more general than the functions studied earlier by Agarwal and Sharma.

\par Finally in 1972, Mittal and Gupta [21] defined a generalized function of two  \linebreak variables popularly known as the H-function of two variables, which generalizes almost all the hypergeometric functions involving  two variables mentioned above.

\par Very recently, the hypergeometric function of two variables introduced by Agarwal [1] and Sharma[27,28] have found interesting applications in \textbf{wireless communications}. For this in a paper by Xia, et.al [32]  by considering dual-hop channel state information (CSI)-assisted amplify-and-forward (AF) relaying over Nakagami-m fading channels, the \linebreak  cumulative distribution function (CDF) of the end-to-end signal-to-noise ratio(SNR) is \linebreak derived . In particular, when the fading shape factors $m_1$ and $m_2$ at consecutive hops take non-integer values, the bivariate H-function and G-function are exploited to obtain an exact analytical expression for the CDF.

\par The remainder of the paper will be organised as follows. In section 2, we shall define the generalization of I-function namely, I-function of two variables. In section 3, we give the notations and results used throughout this paper. In section 4, convergence conditions for this function have been derived. In section 5, we will obtain the series representations and behaviour of the function for small values of  the variables. In section 6, we list special cases of our function giving relations with other functions available in the literature, including H-functions of two variables and G-functions of two variables.  In section 7, we mention few important properties.

\section{The I-function of two variables}

\hspace{20pt} The double Mellin Barnes type contour integral occuring in this paper will be referred to as the I-function of two variables throughout our present study and will be defined and represented in the following manner.\\
\small
 $$ \mathrm{I[z_1,z_2]}= \mathrm{I_{\:p_1,\:q_1 :\; p_2,\:q_2 ;\; p_3,\:q_3}^{\:0,\:n_1:\; m_2,\:n_2;\;m_3,\:n_3}} \left[\begin{array}{c}z_1\\z_2 \end{array} \left|\begin{array}{l} (a_j;\alpha_j,A_j;\xi_j)_{1,p_1} :(c_j,C_j;U_j)_{1,p_2};\; (e_j,E_j;P_j)_{1,p_3}\\ (b_j; \beta_j, B_j; \eta_j)_{1,q_1} :(d_j,D_j ;V_j)_{1,q_2} ;\; (f_j,F_j;Q_j)_{1,q_3} \end{array} \right. \right] $$
$$ = \frac{1}{(2\pi i)^2} \int_{\mathcal{L}_s} \int_{\mathcal{L}_t} \phi(s,t)\; \theta_1(s)\;\theta_2(t)\; z_1^s \;z_2^t \;ds\; dt \qquad \qquad \qquad  \eqno{(2.1)}$$
\normalsize
where $ \phi(s,t) $, $ \theta_1(s) $, $ \theta_2(t) $ are given by  
\small

$$ \phi(s,t) = \frac{\displaystyle \prod_{j=1}^{n_1} \Gamma^{\xi_j}\left(1-a_j+\alpha_js+A_jt\right)}{\displaystyle \prod_{j=n_1 +1}^{p_1} \Gamma^{\xi_j}\left(a_j-\alpha_js-A_jt\right) \; \displaystyle \prod_{j=1}^{q_1}\Gamma^{\eta_j}\left(1-b_j+\beta_js+B_jt\right)} \eqno{(2.2)}$$

$$ \theta_1(s)= \frac{\displaystyle \prod_{j=1}^{n_2} \Gamma^{U_j}\left(1-c_j+C_js\right)\; \displaystyle \prod_{j=1}^{m_2}\Gamma^{V_j}\left(d_j-D_js\right)} {\displaystyle \prod_{j=n_2+1}^{p_2} \Gamma^{U_j}\left(c_j-C_js\right)\; \displaystyle \prod_{j=m_2+1}^{q_2}\Gamma^{V_j}\left(1-d_j+D_js\right)} \eqno{(2.3)}$$
										
$$ \theta_2(t)= \frac{\displaystyle  \prod_{j=1}^{n_3} \Gamma^{P_j}\left(1-e_j+E_jt\right)\; \displaystyle \prod_{j=1}^{m_3}\Gamma^{Q_j}\left(f_j-F_jt\right)} {\displaystyle \prod_{j=n_3+1}^{p_3} \Gamma^{P_j}\left(e_j-E_jt\right)\; \displaystyle \prod_{j=m_3+1}^{q_3}\Gamma^{Q_j}\left(1-f_j+F_jt\right)} \eqno{(2.4)}$$

\normalsize										
Also : 
\begin{itemize}	
\item $ z_1 \neq 0 $ , $ z_2 \neq 0$ ;
\item $i=\sqrt{-1}$; 
\item an empty product is interpreted as unity ;
\item the parameters  $n_j,p_j,q_j(j=1,2,3)$ , $m_j(j=2,3)$  are nonnegative integers  such that $0\leq n_j\leq p_j(j=1,2,3)$ , $q_1\geq 0$ , $0 \leq m_j\leq q_j(j=2,3)$  (not all zero simulataneously) ;
\item  $\alpha_j, A_j (j=1,\ldots,p_1) $ ,  $ \beta_j, B_j (j=1,\ldots,q_1)$ , $ C_j (j=1,\ldots,p_2)$ , \linebreak $D_j (j=1,\ldots,q_2)$ , $E_j (j=1,\ldots,p_3)$  , $F_j (j=1,\ldots,q_3)$ are assumed to be positive quantities for standardisation purpose.

\par However, the definition of I-function of two variables will have a meaning even if some of the quantities are zero or negative numbers. For these we may obtain\linebreak  corresponding transformation formulas which will be given in later section.

\item  $a_j(j=1,\ldots, p_1)$ , $b_j(j=1,\ldots,q_1)$ , $c_j (j=1,\ldots,p_2)$ , $d_j (j=1,\ldots,q_2)$ , $e_j (j=1,\ldots,p_3)$ and  $f_j (j=1,\ldots,q_3)$  are complex numbers;

\item The exponents   $\xi_j (j=1,\ldots,p)$ , $\eta_j (j=1,\ldots,q)$ , $U_j (j=1,\ldots,p_2)$ ,\linebreak  $V_j (j=1,\ldots,q_2)$ , $P_j (j=1,\ldots,p_3)$ , $Q_j (j=1,\ldots,q_3)$ of  various gamma \linebreak functions involved in (2.2) , (2.3) and (2.4)  may take non-integer values.

\item $\mathcal{L}_s$ and  $\mathcal{L}_t$ are suitable contours of Mellin - Barnes type. Moreover, the contour $\mathcal{L}_s$ is in the complex s-plane and runs from $\sigma_1-i\infty$ to $\sigma_1 + i\infty$, ($\sigma_1$ real) so that all the singularities of   $ \Gamma^{V_j}\left(d_j-D_js\right) (j=1,\ldots,m_2)$  lie to the right of $\mathcal{L}_s$   and all the singularities of   $\Gamma^{U_j}\left(1-c_j+C_js\right) (j=1,\ldots,n_2)$ ,  $\Gamma^{\xi_j}\left(1-a_j+\alpha_js+A_jt\right)
\linebreak (j=1,\ldots,n_1)$   lie to the left of $\mathcal{L}_s$ ;

\item  The contour $\mathcal{L}_t$ is in the complex  t-plane and runs from  $\sigma_2-i\infty$  to  $\sigma_2+i\infty$ , ($\sigma_2$ real)  so that all the singularities of  $\Gamma^{Q_j}\left(f_j-F_jt\right) (j=1,\ldots,m_3)$ lie to the right of $\mathcal{L}_t$ , and all the singularities of   $\Gamma^{P_j}\left(1-e_j+E_jt\right)$ $(j=1,\ldots,n_3)$ ,  $\Gamma^{\xi_j}\left(1-a_j+\alpha_js+A_jt\right)$  $(j=1,\ldots,n_1)$ lie to the left of  $\mathcal{L}_t$. 
\end{itemize}
\section{Notations and Results used :}

\hspace{20pt}All the assumptions made above will be retained throughout this paper. For the sake of brevity, the function defined in (2.1) will be simply denoted by either  $\mathrm{I[z_1, z_2]}$  or  $ \mathrm{I\left[\begin{array}{c} z_1\\z_2 \end{array}\right]}$.\\

\par Throughout the present work we shall use the following notations. 
\begin{itemize}
	\item  $(a_j;\alpha_j,A_j;\xi_j)_{1,p}$ stands for     $(a_1;\alpha_1,A_1;\xi_1) , (a_2;\alpha_2,A_2;\xi_2),\ldots,(a_p;\alpha_p,A_p;\xi_p).$
   \item $(c_j,C_j;U_j)_{1,p_2}$ stands for  $(c_1,C_1;U_1), (c_2,C_2;U_2),\ldots,(c_{p_2},C_{p_2};U_{p_2}).$
  \item $(a_j;\alpha_j,A_j;1)_{1,p}$ stands  for $(a_1;\alpha_1,A_1;1) , (a_2;\alpha_2,A_2;1),\ldots, (a_p;\alpha_p,A_p;1).$
\item $(a_j;\alpha_j,A_j)_{1,p}$ stands  for $(a_1;\alpha_1,A_1), (a_2;\alpha_2,A_2),\ldots,(a_p;\alpha_p,A_p).$		
\item $(a_j;\alpha_j,1)_{1,p}$ stands  for $(a_1;\alpha_1,1), (a_2;\alpha_2,1),\ldots,(a_p;\alpha_p,1).$
\item $(a_j;\alpha_j)_{1,p}$ stands  for $(a_1;\alpha_1), (a_2;\alpha_2),\ldots, (a_p;\alpha_p).$
\item $(a_j;1)_{1,p}$ stands  for $(a_1;1), (a_2;1),\ldots,(a_p;1).$
\item $(a_p)=(a_j)_{1,p}$ stands  for $(a_1), (a_2),\ldots, (a_p).$ 
\end{itemize}
 

\section{Convergence Conditions}

\hspace{20pt} Following the results  of Braaksma [6 ,p.278] and Rathie [24], it can easily be shown that the function 
defined by (2.1) is an analytic function of $z_1$ and $z_2$ if $ R<0 $ and $S<0$ where  
\small
$$ R = \quad \sum_{j=1}^{p_1}\xi_j\alpha_j + \sum_{j=1}^{p_2} U_jC_j - \sum_{j=1}^{q_1}\eta_j\beta_j - \sum_{j=1}^{q_2} V_jD_j  \eqno{(4.1)}  $$
$$ S = \quad \sum_{j=1}^{p_1}\xi_jA_j + \sum_{j=1}^{p_3} P_jE_j - \sum_{j=1}^{q_1}\eta_jB_j - \sum_{j=1}^{q_3} Q_jF_j  \eqno{(4.2)} $$
\normalsize
The sufficient conditions for the convergence of (2.1) are given in the following theorem.\\
\\
\textbf{Theorem :}  The integral (2.1) is  convergent   if 
   $|arg(z_1)|< \frac{1}{2}\Delta_1\pi$, $|arg(z_2)|< \frac{1}{2}\Delta_2\pi$ \\ \\and $\Delta_1,\Delta_2 >0$ where 
\small
$$ \Delta_1 = \left[\sum_{j=1}^{n_1}\xi_j\alpha_j - \sum_{j=n_1+1}^{p_1}\xi_j\alpha_j - \sum_{j=1}^{q_1}\eta_j\beta_j + 
\sum_{j=1}^{n_2}U_jC_j \right.\qquad \qquad $$
$$ \qquad \qquad   \left. - \sum_{j=n_2+1}^{p_2}U_jC_j + \sum_{j=1}^{m_2}V_jD_j - \sum_{j=m_2+1}^{q_2}V_jD_j\right] \eqno{(4.3)}$$ 
$$ \Delta_2 = \left[\sum_{j=1}^{n_1}\xi_jA_j - \sum_{j=n_1+1}^{p_1}\xi_jA_j - \sum_{j=1}^{q_1}\eta_jB_j + 
\sum_{j=1}^{n_3}P_jE_j \right.\qquad \qquad $$
$$ \qquad \qquad    \left.- \sum_{j=n_3+1}^{p_3}P_jE_j + \sum_{j=1}^{m_3}Q_jF_j - \sum_{j=m_3+1}^{q_3}Q_jF_j\right] \eqno{(4.4)}$$
\normalsize
and if  $|arg(z_1)|= \frac{1}{2}\Delta_1\pi$ , $|arg(z_2)|= \frac{1}{2}\Delta_2\pi$ and $\Delta_1,\Delta_2 \geq 0$, 
  then integral (2.1)\linebreak  converges absolutely under the following conditions.\\ \\ \textbf{(i)} $\; \mu_1 =\mu_2=0 $, $\Omega_1>1$ and $\Omega_2>1$,where 
	\small 
  $$\mu_1= \sum_{j=1}^{p_2}U_jC_j - \sum_{j=1}^{q_2}V_jD_j + \sum_{j=1}^{p_1} \xi_j\alpha_j - \sum_{j=1}^{q_1}\eta_j\beta_j  \eqno{(4.5)}$$
  $$\mu_2= \sum_{j=1}^{p_3}P_jE_j - \sum_{j=1}^{q_3}Q_jF_j + \sum_{j=1}^{p_1} \xi_jA_j - \sum_{j=1}^{q_1}\eta_jB_j   \eqno{(4.6)}$$
  $$ \Omega_1= \left[\sum_{j=1}^{p_1} \xi_j \left(\Re(a_j)-\frac{1}{2}\right) - \sum_{j=1}^{q_1} \eta_j\left(\Re(b_j)-\frac{1}{2}\right)  \right. \qquad \qquad \qquad  $$
	$$ \qquad    \left. + \sum_{j=1}^{p_2} U_j\left(\Re(c_j)-\frac{1}{2}\right) -\sum_{j=1}^{q_2} V_j\left(\Re(d_j) -\frac{1}{2}\right) \right] \eqno{(4.7)}$$
  $$ \Omega_2= \left[\sum_{j=1}^{p_1} \xi_j\left(\Re(a_j)-\frac{1}{2}\right)- \sum_{j=1}^{q_1} \eta_j\left(\Re(b_j)-\frac{1}{2}\right) \right. \qquad \qquad \qquad $$
	$$ \qquad    \left. +  \sum_{j=1}^{p_3} P_j\left(\Re(e_j)-\frac{1}{2}\right) -\sum_{j=1}^{q_3} Q_j\left(\Re(f_j)-\frac{1}{2}\right)\right] \eqno{(4.8)}  $$  \\
	\normalsize
 \textbf{(ii)}  $\;\mu_1 \neq 0 ,\mu_2\neq 0 $ , if with $ s = \sigma_1+ it_1$ , $t=\sigma_2+it_2$ ($\sigma_1$ , $\sigma_2$ , $t_1$ , $t_2$ are real) , $\sigma_1$ and $\sigma_2$ are chosen so that for $|t_1|\rightarrow\infty$ , $|t_2|\rightarrow\infty$  we have $(\Omega_1+\sigma_1\mu_1) > 1$ and $(\Omega_2+\sigma_2\mu_2) > 1$.\\
   \\
\textbf{Proof:} The existence of  $\mathrm{I[z_1,z_2]}$ may be recognised by the convergence of the integral (2.1) which depends on the asymptotic behaviour of the functions $\phi(s,t)$ , $\theta_1(s)$ and $\theta_2(t)$ defined by (2.2) , (2.3) and (2.4) respectively. Such an asymptotic is based on the following relation for gamma function [18 , 25(p.33)]
\small
$$ |\Gamma(x+iy)| \sim   \sqrt{2\pi} \; |y|^{x-\frac{1}{2}}\;exp\left(-\frac{1}{2}\pi|y|\right), \quad (|y|\rightarrow \infty) \eqno{(4.9)}$$
\normalsize
where x and y are real numbers.  
If we apply  (4.9) for various gamma functions of the integrand of (2.1) with $s = \sigma_1+ it_1$ , $t=\sigma_2+it_2$ where $\sigma_1$ , $\sigma_2$ , $t_1$ , $t_2$ are real. Thus we obtain, for $|t_1|\rightarrow\infty $ and  $|t_2|\rightarrow\infty$ ,  
$$ | \phi(s,t)\: z_1^s\: z_2^t|  \sim  k_1\:|t_1|^{-\Omega_1-\sigma_1\mu_1}\; exp\left(-t_1 (argz_1)-\frac{\pi}{2}|t_1|\Delta_1\right)$$
and 
$$ |\phi(s,t)\: z_1^s\: z_2^t|  \sim  k_2\:|t_2|^{-\Omega_2-\sigma_2\mu_2}\; exp\left(-t_2 (argz_2)-\frac{\pi}{2}|t_2|\Delta_2\right)$$
where $k_1$ and $k_2$ are independent of $t_1$ and $t_2$. Hence the result follows.\\
\\
\textbf{Remark :} If $V_j=1(j=1,\ldots,m_2)$ and $Q_j=1(j=1,\ldots,m_3)$  in (2.1), then the function will be denoted by 
\small
$$ \mathrm{\overline{\;I\;} \left[\begin{array}{c} z_1 \\ z_2 \end{array}\right]}  = \mathrm{I_{\:p_1,\:q_1 :\; p_2,\:q_2 ;\; p_3,\:q_3}^{\:0,\:n_1:\; m_2,\:n_2;\;m_3,\:n_3}} \left[\begin{array}{c}z_1\\z_2 \end{array} \left|\begin{array}{l} (a_j;\alpha_j,A_j;\xi_j)_{1,p_1} :(c_j,C_j;U_j)_{1,n_2}, (c_j,C_j;U_j)_{n_2+1,p_2}; \\ (b_j; \beta_j, B_j; \eta_j)_{1,q_1} :(d_j,D_j ;1)_{1,m_2} , (d_j,D_j ;V_j)_{m_2+1,q_2} ; \end{array} \right. \right. $$
  $$ \qquad \qquad \qquad \qquad \qquad \qquad \qquad \qquad \qquad \qquad \qquad \qquad \left.\begin{array}{l} (e_j,E_j;P_j)_{1,n_3}, (e_j,E_j;P_j)_{n_3+1,p_3} \\  (f_j,F_j; 1)_{1,m_3}, (f_j,F_j;Q_j)_{m_3+1,q_3}  \end{array} \right] $$
	
$$ = \frac{1} {(2\pi i)^2}  \int_{\mathcal{L}_s} \int_{\mathcal{L}_t} \phi(s,t)\; \overline{\theta}_1(s)\;\overline{\theta}_2(t)\; z_1^s \;z_2^t \;ds\; dt  \qquad \qquad \qquad \qquad  \eqno{(4.10)} $$
\\
\normalsize
where  $\phi(s,t)$ is given by (2.2) and $\overline{\theta}_1(s) , \;\overline{\theta}_2(t)$ are given by
\small
$$ \overline{\theta}_1(s)= \frac{\displaystyle \prod_{j=1}^{n_2} \Gamma^{U_j}\left(1-c_j+C_js\right)\; \displaystyle \prod_{j=1}^{m_2}\Gamma \left(d_j-D_js\right)} {\displaystyle \prod_{j=n_2+1}^{p_2} \Gamma^{U_j}\left(c_j-C_js\right)\; \displaystyle \prod_{j=m_2+1}^{q_2}\Gamma^{V_j}\left(1-d_j+D_js\right)} \eqno{(4.11)}$$
										
$$ \overline{\theta}_2(t)= \frac{\displaystyle  \prod_{j=1}^{n_3} \Gamma^{P_j}\left(1-e_j+E_jt\right)\; \displaystyle \prod_{j=1}^{m_3}\Gamma \left(f_j-F_jt\right)} {\displaystyle \prod_{j=n_3+1}^{p_3} \Gamma^{P_j}\left(e_j-E_jt\right)\; \displaystyle \prod_{j=m_3+1}^{q_3}\Gamma^{Q_j}\left(1-f_j+F_jt\right)} \eqno{(4.12)}$$
\normalsize
\section{Series Representation}
\hspace{20 pt} Following the lines of Rathie[24], we can obtain the series representation and \linebreak behaviour for small values for the function  $\mathrm{\overline{I} [z_1 
, z_2]}$ defined and represented by (4.10 ). The series  representation  may be  given as follows :\\
\\If \hspace{10 pt}  \textbf{(i)} $ D_h(d_j + r) \neq D_j(d_h + \mu)$,  for $j \neq h; \; j,h = 1,\ldots,m_2;\; r, \mu = 0,1,2,\ldots $\\ 
 \par \textbf{(ii)}  $ F_l(f_j + k) \neq F_j(f_l + \nu)$,  for $ j \neq l; \; j,l = 1,\ldots,m_3;\; k,\nu = 0,1,2,\ldots$ \\
\par \textbf{(iii)} $z_1 \neq 0 $ ,  $z_2 \neq 0$ , $ R  < 0$ , $S < 0$ (where R  and  S are defined by (4.1) and (4.2) respectively), \\
\par \textbf{(iv)} and if all the poles of (4.10) are simple ,\\
\\
 then the integral (4.10) can be evaluated with the help of the Residue theorem to give \\
\\
\small
$ \mathrm{\overline{\;I\;} [z_1 , z_2]} $\\
$$ = \displaystyle \sum_{h=1}^{m_2} \displaystyle \sum_{l=1}^{m_3} \displaystyle \sum_{r=0}^{\infty} \displaystyle \sum_{k=0}^{\infty} \left\{ \frac{(-1)^{r+k} } { D_h\: F_l \:r!\;  k!} \; \phi\left(\frac{d_h+r}{D_h}, \frac{f_l+k}{F_l}\right) \overline{\theta}_3\left(\frac{d_h+r}{D_h}\right)\;\overline{\theta}_4\left(\frac{f_l+k}{F_l}\right) \; z_1^{\frac{d_h+r}{D_h}} \;  z_2^{\frac{f_l+k}{F_l}}\right\}  \eqno{(5.1)}$$ 
\normalsize
where  $ \phi\left(\frac{d_h+r}{D_h}, \frac{f_l+k}{F_l}\right)$ is defined analogous to  $ \phi(s,t) $  given by (2.2)  and  $ \overline{\theta}_3\left(\frac{d_h+r}{D_h}\right) $ , $ \overline{\theta}_4\left(\frac{f_l+k}{F_l}\right)$ are respectively given by
\small 
$$ \overline{\theta}_3\left(\frac{d_h+r}{D_h}\right) =  \frac{\displaystyle \prod_{j=1}^{n_2} \Gamma^{U_j}\left(1-c_j+C_j\left(\frac{d_h+r}{D_h}\right)\right) \displaystyle \prod_{j=1 \atop (j\neq h)}^{m_2} \Gamma\left(d_j-D_j\left(\frac{d_h+r}{D_h}\right)\right)} {\displaystyle \prod_{j=n_2+1}^{p_2} \Gamma^{U_j}\left(c_j - C_j\left(\frac{d_h+r}{D_h}\right)\right) \displaystyle \prod_{j= m_2+1}^{q_2} \Gamma^{V_j}\left(1-d_j+D_j\left(\frac{d_h+r}{D_h}\right)\right)} \eqno{(5.2)}$$ 
$$ \overline{\theta}_4\left(\frac{f_l+k}{F_l}\right)= \frac{\displaystyle \prod_{j=1}^{n_3} \Gamma^{P_j}\left(1-e_j+E_j\left(\frac{f_l+k}{F_l}\right) \right) \displaystyle \prod_{j=1 \atop (j\neq l)}^{m_3} \Gamma\left(f_j-F_j\left(\frac{f_l+k}{F_l}\right)\right)} {\displaystyle \prod_{j=n_3+1}^{p_3} \Gamma^{P_j}\left(e_j-E_j\left(\frac{f_l+k}{F_l}\right)\right) \displaystyle \prod_{j=m_3+1}^{q_3} \Gamma^{Q_j}\left(1-f_j+F_j\left(\frac{f_l+k}{F_l}\right)\right)}     \eqno{(5.3)} $$
\normalsize
for $|z_1| < 1\; , \; |z_2| < 1 $.\\
\\
From (5.1) it follows that \\
$$\mathrm{\overline{\;I\;} [z_1 , z_2]} =  \emph{O} (|z_1|^\alpha \; |z_2|^\beta) ,  \qquad  max{\left\{|z_1|,\; |z_2|\right\}} \rightarrow 0  \eqno{(5.4)} $$
where $$ \alpha\; = \quad  \displaystyle {min \atop {1\leq j\leq m_2}}  \left[\Re\left(\frac{d_j}{D_j}\right)\right]   \eqno{(5.5)} $$
 $$ \beta\; =   \quad \displaystyle { min \atop { 1\leq j\leq m_3}} \left[\Re\left(\frac{f_j}{F_j}\right)\right]   \eqno{(5.6)} $$
for small values of $z_1$ and $z_2$.
\section{Elementary Special Cases}
\hspace{20pt} In this section, we mention some interesting and useful special cases of the I-function of two variables. \\\\
\textbf{(i)} If all the exponents  $\xi_j (j=1,\ldots,p_1)$ , $\eta_j (j=1,\ldots,q_1)$ , $U_j (j=1,\ldots,p_2)$ , $V_j (j=1,\ldots,q_2)$,  $ P_j (j=1,\ldots,p_3)$ and  $ Q_j (j=1,\ldots,q_3) $  are equal to unity , then (2.1)  reduces to the H-function of two variables defined by Mittal and Gupta [21].\\\\
\textbf{(ii)} If we take $p_1=q_1=n_1=0$ in (2.1) then it degenerates into the product of two \linebreak I- functions of one variable introduced by Rathie [24] as \\ \\
\small
$ \quad \mathrm{{I}_{\:0,\:0 :\; p_2,\:q_2;\; p_3,\:q_3}^{\:0,\:0:\; m_2,\:n_2;\;m_3,\:n_3}} \left[\begin{array}{c}z_1\\z_2 \end{array} \left|\begin{array}{l} \textbf{-----} :(c_j,C_j;U_j)_{1,p_2} ;\; (e_j,E_j;P_j)_{1,p_3}\\ \textbf{-----} :(d_j,D_j ;V_j)_{1,q_2} ;\; (f_j,F_j;Q_j)_{1,q_3} \end{array} \right. \right] \qquad \qquad  $ 
$$ \qquad   = \mathrm{{I}_{\:p_2,\:q_2}^{\:m_2,\:n_2}} \left[\begin{array}{c}z_1 \end{array} \left|\begin{array}{l}  (c_j,C_j;U_j)_{1,p_2}\\ (d_j,D_j ;V_j)_{1,q_2} \end{array} \right. \right] \quad \times \quad  \mathrm{I_{\: p_3,\:q_3}^{\:m_3,\:n_3}} \left[\begin{array}{c} z_2 \end{array} \left|\begin{array}{l} (e_j,E_j;P_j)_{1,p_3} \\ (f_j,F_j;Q_j)_{1,q_3} \end{array} \right. \right] \eqno{(6.1)}  $$ 
 \normalsize
\textbf{(iii)} If we take $m_3=1$ , $n_3=p_3$ , $f_1=0$ , $(F_j)_{1,q_3}=1$ , $(A_j)_{1,p_1}=(B_j)_{1,q_1}=(E_j)_{1,p_3}=(F_j)_{1,q_3}=1$ , equate the exponents $P_j(j=1,\ldots,p_3), Q_j(j=1,\ldots,q_3)$ to unity,   replace $q_3$ by $q_3+1$  and let $z_2\rightarrow 0$ in (2.1), we get the following relation by virtue of (2.1) and  known results[10 , p.208],\\
\\
\small
$ \stackrel{lim}{{z_2\rightarrow 0}} \mathrm{{I}_{\:p_1,\:q_1 :\; p_2,\:q_2;\; p_3,\:q_3}^{\:0,\:n_1:\; m_2,\:n_2;\;m_3,\:n_3}} \left[\begin{array}{c}z_1\\z_2\end{array} \left|\begin{array}{l} (a_j;\alpha_j,1;\xi_j)_{1,p_1} : (c_j,C_j;U_j)_{1,p_2}; (e_j,1;1)_{1,p_3}\\ (b_j; \beta_j, 1; \eta_j)_{1,q_1} :(d_j,D_j ;V_j)_{1,q_2};(0,1;1), (f_j,1;1)_{2,q_3+1} \end{array} \right. \right] $
$$ = \frac{\prod_{j=1}^{p_3}\Gamma(1-e_j)}{\prod_{j=1}^{q_3}\Gamma(1-f_j)}\; 
   \mathrm{I_{\:p_1+p_2,\;q_1+q_2}^{\:m_2,n_1+n_2}} \left[\begin{array}{c}z_1 \end{array} \left|\begin{array}{l} (a_j,\alpha_j;\xi_j)_{1,n_1},(c_j,C_j;U_j)_{1,p_2}, (a_j,\alpha_j;\xi_j)_{n_1+1,p_1}\\ (d_j,D_j ;V_j)_{1,q_2}, (b_j, \beta_j; \eta_j)_{1,q_1}  \end{array} \right. \right] \eqno{(6.2)} $$
\normalsize
where $p_1+p_3 < q_1+q_3+1$ \\ \\
\textbf{(iv)} 
A relationship between the I-function of two variables  and the G-function of two \linebreak variables [30, p.7, eq.1.2.3]  is, \\\\
\small 
 $ \mathrm{{I}_{\:p_1,\:q_1 :\; p_2,\:q_2;\; p_3,\:q_3}^{\:0,\:n_1:\; m_2,\:n_2;\;m_3,\:n_3}} \left[\begin{array}{c}z_1\\z_2 \end{array} \left|\begin{array}{l} (a_j;1,1;1)_{1,p_1} :(c_j,1;1)_{1,p_2};\; (e_j,1;1)_{1,p_3}\\ (b_j; 1,1 ; 1)_{1,q_1} :(d_j,1 ; 1)_{1,q_2} ; \;(f_j, 1 ;1)_{1,q_3} \end{array} \right. \right] $
 $$ = \mathrm {G\:(z_1,z_2)} = \mathrm{{G}_{\:p_1,\:q_1 :\; p_2,\:q_2;\; p_3,\:q_3}^{\:0,\:n_1:\; m_2,\:n_2;\;m_3,\:n_3}}  \left[\begin{array}{c}z_1\\z_2 \end{array} \left|\begin{array}{l} (a_{p_1}) :(c_{p_2});\; (e_{p_3})\\ (b_{q_1}) :(d_{q_2}) ;\; (f_{q_3}) \end{array} \right. \right] \qquad \qquad \qquad \eqno{(6.3)}$$ 
\normalsize
\\
\textbf{(v)} Another specialization of parameters in the I-function of two variables yields the \linebreak interesting relationship : \\\\
\small
 $ \mathrm{{I}_{\:p_1,\:q_1 :\; p_2,\:q_2+1;\; p_3,\:q_3+1}^{\:0,\:p_1:\; 1,\:p_2;\;1,\:p_3}} \left[\begin{array}{c} -z_1\\ -z_2 \end{array} \left|\begin{array}{l} (1-a_j;\alpha_j,A_j;1)_{1,p_1} :\\ (1-b_j; \beta_j, B_j; 1)_{1,q_1} : \end{array} \right. \right. $   
        $$ \qquad \qquad \qquad \qquad \qquad \qquad  \left.\begin{array}{l} (1-c_j,C_j; 1)_{1,p_2}; (1-e_j,E_j; 1)_{1,p_3}\\ (0,1;1),(1-d_j,D_j ; 1)_{1,q_2} ; (0,1;1),(1-f_j,F_j; 1)_{1,q_3} \end{array} \right] $$ 				
$$= \mathrm{S_{\:q_1;\:q_2;\:q_3}^{\:p_1;\:p_2;\:p_3}}  \left[ \begin{array}{l} (a_j;\alpha_j,A_j)_{1,p_1} :(c_j,C_j)_{1,p_2}; (e_j,E_j)_{1,p_3} \; ;\\ (b_j; \beta_j, B_j)_{1,q_1} :(d_j,D_j )_{1,q_2} ; (f_j,F_j)_{1,q_3} \; ; \end{array}  \begin{array}{c} z_1 ,  z_2  \end{array} \right] \eqno{(6.4)}$$ 
 where $ \mathrm{S (z_1, z_2)} $ is   generalized Kampe de Feriet function  defined by Srivastava and Daoust[29, 30.p.6,eq.1.2.2].\\\\
\normalsize
\textbf{(vi)}
\small
$ \mathrm{{I}_{\:p_1,\:q_1 :\; p_2,\:q_2+1;\; p_3,\:q_3+1}^{\:0,\:n_1:\; 1,\:p_2;\;1,\:p_3}} \left[\begin{array}{c} z_1\\ z_2 \end{array} \left|\begin{array}{l} (1-a_j;1,1;1)_{1,p_1} :\\ (1-b_j; 1, 1; 1)_{1,q_1} : \end{array} \right. \right. \qquad \qquad \qquad \qquad $ \\ 
$$ \qquad \qquad \qquad \qquad \qquad \qquad \qquad  \left. \begin{array}{l}(1-c_j,1; 1)_{1,p_2};  (1- e_j,1; 1)_{1,p_3}\\ (0,1;1),(1-d_j, 1 ; 1)_{1,q_2} ; (0,1;1), (1-f_j, 1; 1)_{1,q_3}\end{array} \right] $$
 $$  \qquad     = \frac{\prod_{j=1}^{p_1} \Gamma\left(a_j\right) \;  \prod_{j=1}^{p_2} \Gamma\left(c_j\right) \;  \prod_{j=1}^{p_3} \Gamma\left(e_j\right)} {\prod_{j=1}^{q_1} \Gamma \left(b_j\right) \;  \prod_{j=1}^{q_2} \Gamma \left(d_j\right) \;  \prod_{j=1}^{q_3} \Gamma \left(f_j\right)}  \quad \times  \qquad \qquad \qquad \qquad \qquad \qquad \qquad  \qquad \qquad \qquad \qquad $$
$$    \mathrm{F_{\:q_1: q_2; q_3}^{\:p_1 : p_2; p_3}} \left[\begin{array}{c}  (a_j)_{1,p_1} :(c_j)_{1,p_2}, (e_j)_{1,p_3} \quad ; z_1\\ (b_j)_{1,q_1} : (d_j)_{1,q_2} ; (f_j)_{1,q_3} \quad ; z_2  \end{array}  \right] \qquad  \eqno{(6.5)}$$ 
\normalsize  
where the function $\mathrm{F[z_1,z_2]}$ is the Kamp$\acute{e}$ de F$\acute{e}$riet function [30,  p.5, eq.1.2.1]. \\
\\
\textbf{(vii)}
\small
$\mathrm{{I}_{\:0,\:0 :\; p_2,\:q_2+1;\; p_3,\:q_3+1}^{\:0,\:0:\; 1,\:p_2;\; 1,\:p_3}} \left[\begin{array}{c} -z_1\\ -z_2 \end{array} \left|\begin{array}{l} \textbf{----------} :\\ \textbf{----------} : \end{array} \right. \right.  $  
$$ \qquad \qquad \qquad \qquad \qquad \qquad \qquad  \left.\begin{array}{l} (1-c_j,C_j;1)_{1,p_2}; (1-e_j,E_j; 1)_{1,p_3} \\ (0,1;1),(1-d_j,D_j ; 1)_{1,q_2} ; (0,1;1), (1-f_j,F_j;1)_{1,q_3}  \end{array} \right] $$
$$ = \;_{p_2}\Psi_{q_2}   \left[\begin{array}{c} (c_j, C_j)_{1,p_2} \\ (d_j, D_j )_{1,q_2} \end{array} \begin{array} {c} ; \;z_1 \end{array} \right] \quad  \times  \quad _{p_3}\Psi_{q_3}   \left[\begin{array}{c} (e_j, E_j)_{1,p_3} \\ (f_j, F_j )_{1,q_3} \end{array} \begin{array} {c} ; \; z_2 \end{array} \right] \eqno{(6.6)}$$		
\normalsize
 where the function $_{p_2}\Psi_{q_2} $  and  $_{p_3}\Psi_{q_3} $  are   Wright's generalized  hypergeometric \linebreak  functions.[30, p. 19,eq.2.6.11]\\
\\
\textbf{(viii)} 
\small
 $\mathrm{{I}_{\:0,\;0 :\;0,\:2;\; 0,\: 2}^{\:0,\:0:\; 1,\:0;\;1,0}} \left[\begin{array}{c} z_1\\ z_2 \end{array} \left|\begin{array}{l} \textbf{-----}\quad  : \qquad \qquad  \textbf{-----}\quad ; \qquad \qquad  \textbf{-----}\\ \textbf{-----}\quad  :(0,1;1),  (-\mu, \alpha ; 1) ; (0,1;1), (-\nu,\beta; 1) \end{array} \right. \right] $\\   
$$ = \quad J_\mu^{\alpha}(z_1) \quad \times \quad  J_\nu^{\beta} (z_2) \qquad \qquad \qquad \qquad \qquad   \eqno{(6.7)}$$
\normalsize
where the functions $ J_\mu^\alpha(z_1)$ and  $J_\nu^\beta (z_2)$ are  Wright's generalized Bessel functions.[30, p.19,eq.2.6.10] \\
\\
\textbf{(ix)}
\small
 $ \mathrm{{I}_{\:0,\:0 :\; 2,\:2;\;2,\: 2}^{\:0,\:0:\; 1,\:2;\;1,\:2}} \left[\begin{array}{c} -z_1\\ -z_2 \end{array} \left|\begin{array}{l} \textbf{---} : (1,1;1),(1-\alpha,1;p); (1,1;1),(1-\beta,1;q)\\ \textbf{---} :(0,1;1),  (-\alpha, 1; p) ; (0,1;1), (-\beta,1; q) \end{array} \right. \right]  $ \\
 $$  = \quad \Phi(z_1,p,\alpha) \quad \times \quad \Phi(z_2,q,\beta) \qquad \qquad \qquad \eqno{(6.8)}$$
\normalsize
where $\Phi(z_1,p,\alpha)$  and   $\Phi(z_2,q,\beta)$ are the generalised Riemann-zeta functions[10, p.27, 1.11, eq(1)] , which are generalizations of Hurwitz zeta functions and Riemann zeta \linebreak  functions [10, p.24, 1.10, eq(1) and 1.12, eq(1)] \\\\
\textbf{(x)}
\small
$
\mathrm{{I}_{\:0,\:0 :\; 2,\:2;\;2,\: 2}^{\:0,\:0:\; 1,\:2;\;1,\:2}} \left[\begin{array}{c} -z_1\\ -z_2 \end{array} \left|\begin{array}{l} \textbf{---} : (1,1;1),(1,1;p); (1,1;1),(1,1;q)\\ \textbf{---} :(0,1;1),  (0, 1; p) ; (0,1;1), (0,1; q) \end{array} \right. \right]  $ 
$$ =  \quad  F(z_1,p)\quad \times \quad  F(z_2,q) \qquad \qquad \qquad \qquad \qquad \qquad \qquad  \eqno{(6.9)}  $$
\normalsize
where $F(z_1,p)$  and $F(z_2,q)$ are the polylogarithms of order p and q respectively. For p=2 and q=2 , the R.H.S. of (6.9) equation reduces to product of Eulers's dilogarithm [10, p.31, 1.11.1, eq(2)]
\normalsize
\section{Elementary properties and Transformation formulas}
The properties given below are immediate consequence of the definition (2.1) and hence they are given here without proof:\\\\
\textbf{(i)}\\
\small
$  \mathrm{{I}_{\:p_1,\:q_1 :\; p_2,\:q_2;\; p_3,\:q_3}^{\:0,\:0:\; m_2,\:n_2;\;m_3,\:n_3}} \left[\begin{array}{c}z_1\\z_2 \end{array} \left|\begin{array}{l} (a_j;\alpha_j,A_j;\xi_j)_{1,p_1} :(c_j,C_j;U_j)_{1,p_2}; (e_j,E_j;P_j)_{1,p_3}\\ (b_j; \beta_j, B_j; \eta_j)_{1,q_1} :(d_j,D_j ;V_j)_{1,q_2} ; (f_j,F_j;Q_j)_{1,q_3} \end{array} \right. \right]  \qquad \qquad \qquad \qquad $\\\\
$\quad  = \mathrm{{I}_{\:q_1,\:p_1 :\; q_2,\:p_2;\; q_3,\:p_3}^{\:0,\:0:\; n_2,\:m_2;\;n_3,\:m_3}} \left[\begin{array}{c}z_1^{-1}\\z_2^{-1} \end{array} \left|\begin{array}{l}(1-b_j; \beta_j, B_j; \eta_j)_{1,q_1} : \\ (1-a_j;\alpha_j,A_j;\xi_j)_{1,p_1} : \end{array} \right. \right.$\\ 
 $$\qquad \qquad \qquad \qquad \qquad \qquad \qquad  \left. \begin{array}{l} (1-d_j,D_j ;V_j)_{1,q_2} ; (1-f_j,F_j;Q_j)_{1,q_3}\\ (1-c_j,C_j;U_j)_{1,p_2}; (1-e_j,E_j;P_j)_{1,p_3} \end{array} \right] \eqno{(7.1)}$$
\normalsize
\textbf{(ii)}\\ 
\small
$  z_1^{k_1}\; z_2^{k_2}\;  \mathrm{I [z_1,z_2]} $\\ 
$ = \mathrm{{I}_{\:p_1,\:q_1 :\; p_2,\:q_2;\; p_3,\:q_3}^{\:0,\:n_1:\; m_2,\:n_2;\;m_3,\:n_3}}\left[\begin{array}{l}z_1\\z_2 \end{array} \left|\begin{array}{l} (a_j+k_1\alpha_j+ k_2A_j;\alpha_j, A_j;\xi_j)_{1,p_1} : \\ (b_j+k_1\beta_j+k_2B_j; \beta_j, B_j; \eta_j)_{1,q_1} :   \end{array} \right.  \right.\qquad \qquad $\\ 
$$ \qquad \qquad \qquad  \qquad \qquad \qquad  \left.\begin{array}{l} (c_j+ k_1C_j,C_j;U_j)_{1,p_2}; (e_j+k_2E_j,E_j;P_j)_{1,p_3} \\ (d_j+k_1D_j,D_j ;V_j)_{1,q_2} ; (f_j+k_2F_j,F_j;Q_j)_{1,q_3} \end{array} \right]   \eqno{(7.2)} $$ 
\normalsize
for $k_1 > 0 ,\; k_2 > 0 $.\\\\
\textbf{(iii)}\\
\small 
$ \frac{1}{k_1} \;\frac{1}{k_2} \; \mathrm{I [z_1,z_2]} $\\
$$=  \mathrm{{I}_{\:p_1,\:q_1 :\; p_2,\:q_2;\; p_3,\:q_3}^{\:0,\:n_1:\; m_2,\:n_2;\;m_3,\:n_3}}\left[\begin{array}{c}z_1^{k_1}\\z_2^{k_2} \end{array} \left|\begin{array}{l} (a_j; k_1\alpha_j, k_2A_j;\xi_j)_{1,p_1} :(c_j,k_1C_j;U_j)_{1,p_2}; (e_j, k_2E_j;P_j)_{1,p_3}\\ (b_j;k_1\beta_j,k_2B_j; \eta_j)_{1,q_1} :(d_j, k_1D_j ;V_j)_{1,q_2} ; (f_j,k_2F_j ;Q_j)_{1,q_3} \end{array} \right. \right]\eqno{(7.3)} $$
\normalsize
where $k_1 > 0 ,\; k_2 > 0 $.\\\\
\textbf{(iv)}\\
\small
 $ \mathrm{{I}_{\:p_1,\:q_1 :\; p_2,\:q_2;\; p_3,\:q_3}^{\:0,\:n_1:\; m_2,\:n_2,\;m_3,\:n_3}} \left[\begin{array}{c}z_1\\z_2 \end{array} \left|\begin{array}{l} (a;\alpha,0;\xi), (a_j;\alpha_j,A_j;\xi_j)_{2,p_1} :(c_j,C_j;U_j)_{1,p_2}; (e_j,E_j;P_j)_{1,p_3}\\ (b_j; \beta_j, B_j; \eta_j)_{1,q_1} :(d_j,D_j ;V_j)_{1,q_2} , (f_j,F_j;Q_j)_{1,q_3} \end{array} \right. \right] $ \\\\
$ =  \mathrm{{I}_{\:p_1-1,\:q_1 :\; p_2+1,\:q_2;\; p_3\;q_3}^{\:0,\:n_1-1:\; m_2,\:n_2+1;\;m_3,\:n_3}} \left[\begin{array}{c}z_1\\z_2 \end{array} \left|\begin{array}{l} (a_j;\alpha_j,A_j;\xi_j)_{2,p_1} : \\ (b_j; \beta_j, B_j; \eta_j)_{1,q_1} : \end{array} \right. \right. $\\ 
$$ \qquad \qquad \qquad \qquad \qquad \qquad  \left. \begin{array} {l} (a,\alpha;\xi),(c_j,C_j;U_j)_{1,p_2}; (e_j,E_j;P_j)_{1,p_3}\\  (d_j,D_j ;V_j)_{1,q_2} ; (f_j,F_j;Q_j)_{1,q_3} \end{array} \right] \eqno{(7.4)} $$
\normalsize
where $p_1 \geq n_1 \geq 1$.\\
\\
\textbf{(v)}\\
\small
 $ \mathrm{{I}_{\:p_1,\:q_1 :\; p_2,\:q_2;\; p_3,\:q_3}^{\:0,\:n_1:\; m_2,\:n_2;\;m_3,\:n_3}} \left[\begin{array}{c}z_1\\z_2 \end{array} \left|\begin{array}{l}  (a_j;\alpha_j,A_j;\xi_j)_{1,p_1-1} , (a;\alpha,0;\xi) :(c_j,C_j;U_j)_{1,p_2}; (e_j,E_j;P_j)_{1,p_3}\\ (b_j; \beta_j, B_j; \eta_j)_{1,q_1} :(d_j,D_j ;V_j)_{1,q_2} ; (f_j,F_j;Q_j)_{1,q_3} \end{array} \right. \right] $ \\\\
$ =  \mathrm{{I}_{\:p_1-1,\:q_1 :\; p_2+1,\:q_2;\; p_3,\:q_3}^{\:0,\:n_1:\; m_2,\:n_2,\:m_3,\:n_3}} \left[\begin{array}{c}z_1\\z_2 \end{array} \left|\begin{array}{l} (a_j;\alpha_j,A_j;\xi_j)_{1,p_1-1} :\\ (b_j; \beta_j, B_j; \eta_j)_{1,q_1} : \end{array} \right. \right.$\\
 $$ \qquad \qquad \qquad  \qquad \qquad \qquad \left. \begin{array}{l} (c_j,C_j;U_j)_{1,p_2},(a,\alpha;\xi); (e_j,E_j;P_j)_{1,p_3}\\  (d_j,D_j ;V_j)_{1,q_2} ; (f_j,F_j;Q_j)_{1,q_3} \end{array} \right] \eqno{(7.5)} $$
\normalsize
where  $ p_1-1 \geq n_1 \geq 0$ \\\\
\textbf{(vi)}\\
\small
 $ \mathrm{{I}_{\:p_1,\:q_1 :\; p_2,\:q_2;\; p_3,\:q_3}^{\:0,\:n_1:\; m_2,\:n_2;\;m_3,\:n_3} }\left[\begin{array}{c}z_1\\z_2 \end{array} \left|\begin{array}{l}  (a_j;\alpha_j,A_j;\xi_j)_{1,p_1} :(c_j,C_j;U_j)_{1,p_2}; (e_j,E_j;P_j)_{1,p_3}\\ (b_j; \beta_j, B_j; \eta_j)_{1,q_1-1}, (b; \beta,0;\eta) :(d_j,D_j ;V_j)_{1,q_2} ; (f_j,F_j;Q_j)_{1,q_3} \end{array} \right. \right] $ \\\\
$ =  \mathrm{{I}_{\:p_1,\:q_1-1 :\; p_2,\:q_2+1;\; p_3,\:q_3}^{\:0,\:n_1:\; m_2,\:n_2;\;m_3,\:n_3}} \left[\begin{array}{c}z_1\\z_2 \end{array} \left|\begin{array}{l} (a_j;\alpha_j,A_j;\xi_j)_{1,p_1} :\\ (b_j; \beta_j, B_j; \eta_j)_{1,q_1-1} :  \end{array} \right. \right. $\\ 
 $$\qquad  \qquad \qquad \qquad  \qquad \qquad   \left. \begin{array}{l} (c_j,C_j;U_j)_{1,p_2}; (e_j,E_j;P_j)_{1,p_3}\\ (d_j,D_j ;V_j)_{1,q_2}, (b,\beta;\eta) ; (f_j,F_j;Q_j)_{1,q_3} \end{array} \right] \eqno{(7.6)} $$ 
\normalsize
where $q_1 - 1 \geq 0$ \\\\
\textbf{(vii)}\\
\small
$   \mathrm{{I}_{\:p_1,\:q_1 :\; p_2,\:q_2;\; p_3,\:q_3}^{\:0,\:n_1:\;m_2,\:n_2;\:m_3,\:n_3}} \left[\begin{array}{c}z_1\\z_2 \end{array} \left|\begin{array}{l} (a;0,0;\xi), (a_j;\alpha_j,A_j;\xi_j)_{2,p_1} :(c_j,C_j;U_j)_{1,p_2}; (e_j,E_j;P_j)_{1,p_3}\\ (b_j; \beta_j, B_j; \eta_j)_{1,q_1} :(d_j,D_j ;V_j)_{1,q_2} ; (f_j,F_j;Q_j)_{1,q_3} \end{array} \right. \right] $\\ \\ 
$ = \quad \Gamma^{\xi}\left(1-a\right)  \quad \times $  \\ 
$$\mathrm{{I}_{\:p_1-1,\:q_1 :\; p_2,\:q_2;\; p_3,\:q_3}^{\:0,\:n_1-1:\; m_2,\:n_2;\;m_3,\:n_3}} \left[\begin{array}{c}z_1\\z_2 \end{array} \left|\begin{array}{l} (a_j;\alpha_j,A_j;\xi_j)_{2,p_1} :(c_j,C_j;U_j)_{1,p_2};(e_j,E_j;P_j)_{1,p_3}\\ (b_j; \beta_j, B_j; \eta_j)_{1,q_1} :(d_j,D_j ;V_j)_{1,q_2} ; (f_j,F_j;Q_j)_{1,q_3} \end{array} \right. \right] \eqno{(7.7)}$$
\normalsize
where $ \geq n_1 \geq 1 $ , $ \Re(1-a)>0$. \\\\
\textbf{(viii)}\\
\small
$\mathrm{{I}_{\:p_1,\:q_1 :\; p_2,\:q_2;\; p_3,\:q_3}^{\:0,\:n_1:\; m_2,\:n_2;,\;m_3,\:n_3}} \left[\begin{array}{c}z_1\\z_2 \end{array} \left|\begin{array}{l} (a_j;\alpha_j,A_j;\xi_j)_{1,p_1-1}, (a;0,0;\xi) :(c_j,C_j;U_j)_{1,p_2}; (e_j,E_j;P_j)_{1,p_3}\\ (b_j; \beta_j, B_j; \eta_j)_{1,q_1} :(d_j,D_j ;V_j)_{1,q_2} ; (f_j,F_j;Q_j)_{1,q_3} \end{array}\right. \right] $\\\\ 
$=\quad  \frac{1}{\Gamma^{\xi}\left(a \right)} \quad \times $\\
 $$ \mathrm{{I}_{\:p_1-1,\:q_1 :\; p_2,\:q_2;\; p_3,\:q_3}^{\:0,\:n_1:\; m_2,\:n_2;\;m_3,\:n_3}} \left[\begin{array}{c}z_1\\z_2 \end{array} \left|\begin{array}{l} (a_j;\alpha_j,A_j;\xi_j)_{1,p_1-1} :(c_j,C_j;U_j)_{1,p_2}, (e_j,E_j;P_j)_{1,p_3}\\ (b_j; \beta_j, B_j; \eta_j)_{1,q_1} :(d_j,D_j ;V_j)_{1,q_2} , (f_j,F_j;Q_j)_{1,q_3} \end{array} \right. \right] \eqno{(7.8)}$$
\normalsize
where  $p_1-1 \geq n_1 \geq 0$ , $\Re(a) > 0 $\\\\
\textbf{(ix)}\\
\small
$   \mathrm{{I}_{\:p_1,\:q_1 :\; p_2,\:q_2;\; p_3,\:q_3}^{\:0,\:n_1:\; m_2,\:n_2;\;m_3,\:n_3}} \left[\begin{array}{c}z_1\\z_2 \end{array} \left|\begin{array}{l} (a_j;\alpha_j,A_j;\xi_j)_{1,p_1} :(c_j,C_j;U_j)_{1,p_2}; (e_j,E_j;P_j)_{1,p_3}\\ (b_j; \beta_j, B_j; \eta_j)_{1,q_1-1}, (b;0,0; \eta) :(d_j,D_j ;V_j)_{1,q_2} ; (f_j,F_j;Q_j)_{1,q_3} \end{array} \right. \right] $\\\\ 
$= \quad \frac{1}{\Gamma^{\eta}\left(1-b\right)} \quad \times $ \\
$$\mathrm{{I}_{\:p_1,\:q_1-1 :\; p_2,\:q_2;\; p_3,\:q_3}^{\:0,\:n_1:\; m_2,\:n_2;\;m_3,\:n_3} }\left[\begin{array}{c}z_1\\z_2 \end{array} \left|\begin{array}{l} (a_j;\alpha_j,A_j;\xi_j)_{1,p_1} :(c_j,C_j;U_j)_{1,p_2}; (e_j,E_j;P_j)_{1,p_3}\\ (b_j; \beta_j, B_j; \eta_j)_{1,q_1-1} :(d_j,D_j ;V_j)_{1,q_2} ; (f_j,F_j;Q_j)_{1,q_3} \end{array} \right. \right] \eqno{(7.9)}$$
\normalsize
where  $ q_1 - 1 \geq 0 $, $ \Re(1-b) > 0 $\\\\
\textbf{(x)}\\
\small
$   \mathrm{{I}_{\:p_1,\:q_1 :\; p_2,\:q_2;\; p_3,\:q_3}^{\:0,\:n_1:\; m_2,\:n_2;\;m_3,\:n_3}} \left[\begin{array}{c}z_1\\z_2 \end{array} \left|\begin{array}{l} (a_j;\alpha_j,A_j;\xi_j)_{1,p_1} :(c,0;U),(c_j,C_j;U_j)_{2,p_2}; (e_j,E_j;P_j)_{1,p_3}\\ (b_j; \beta_j, B_j; \eta_j)_{1,q_1} :(d_j,D_j ;V_j)_{1,q_2} ; (f_j,F_j;Q_j)_{1,q_3} \end{array} \right. \right] $\\\\ 
$= \quad \Gamma^{U}\left(1-c\right) \quad \times $ \\
$$ \mathrm{{I}_{\:p_1,\:q_1 :\; p_2-1,\:q_2;\; p_3,\:q_3}^{\:0,\:n_1:\; m_2,\:n_2-1;\;m_3,\:n_3}} \left[\begin{array}{c}z_1\\z_2 \end{array} \left|\begin{array}{l} (a_j;\alpha_j,A_j;\xi_j)_{1,p_1} :(c_j,C_j;U_j)_{2,p_2}; (e_j,E_j;P_j)_{1,p_3}\\ (b_j; \beta_j, B_j; \eta_j)_{1,q_1} :(d_j,D_j ;V_j)_{1,q_2}; (f_j,F_j;Q_j)_{1,q_3} \end{array} \right. \right] \eqno{(7.10)}$$
\normalsize
where  $ p_2 \geq n_2 \geq 1$ , $ \Re(1-c) > 0 $\\\\
\textbf{(xi)}\\
\small
$ \mathrm{{I}_{\:p_1,\:q_1 :\; p_2,\:q_2;\; p_3,\:q_3}^{\:0,\:n_1:\; m_2,\:n_2;\;m_3,\:n_3}} \left[\begin{array}{c}z_1\\z_2 \end{array} \left|\begin{array}{l} (a_j;\alpha_j,A_j;\xi_j)_{1,p_1} :(c_j,C_j;U_j)_{1,p_2-1},(c,0;U); (e_j,E_j;P_j)_{1,p_3}\\ (b_j; \beta_j, B_j; \eta_j)_{1,q_1} :(d_j,D_j ;V_j)_{1,q_2} ; (f_j,F_j;Q_j)_{1,q_3} \end{array} \right. \right] $\\ \\
$= \quad \frac{1} {\Gamma^{U}\left(c\right)} \quad \times $\\
 $$\mathrm{{I}_{\:p_1,\:q_1 :\; p_2-1,\:q_2;\; p_3,\:q_3}^{\:0,\:n_1:\; m_2,\:n_2;\;m_3,\:n_3}} \left[\begin{array}{c}z_1\\z_2 \end{array} \left|\begin{array}{l} (a_j;\alpha_j,A_j;\xi_j)_{1,p_1} :(c_j,C_j;U_j)_{1,p_2-1}; (e_j,E_j;P_j)_{1,p_3}\\ (b_j; \beta_j, B_j; \eta_j)_{1,q_1} :(d_j,D_j ;V_j)_{1,q_2} ; (f_j,F_j;Q_j)_{1,q_3} \end{array} \right. \right] \eqno{(7.11)}$$
\normalsize
where $p_2 - 1 \geq n_2 \geq 0 $ , $ \Re(c) > 0 $\\\\
\textbf{(xii)}\\
\small
$   \mathrm{{I}_{\:p_1,\:q_1 :\; p_2,\:q_2;\; p_3,\:q_3}^{\:0,\:n_1:\; m_2,\:n_2;\;m_3,\:n_3}} \left[\begin{array}{c}z_1\\z_2 \end{array} \left|\begin{array}{l} (a_j;\alpha_j,A_j;\xi_j)_{1,p_1} :(c_j,C_j;U_j)_{1,p_2}; (e_j,E_j;P_j)_{1,p_3}\\ (b_j; \beta_j, B_j; \eta_j)_{1,q_1} :(d,0;V), (d_j,D_j ;V_j)_{2,q_2} ; (f_j,F_j;Q_j)_{1,q_3} \end{array} \right. \right] $ \\ \\
$= \quad \Gamma^{V}\left(d\right) \quad \times $\\
$$\mathrm{{I}_{\:p_1,\:q_1 :\; p_2,\:q_2-1;\; p_3,\:q_3}^{\:0,\:n_1:\; m_2-1,\:n_2;\;m_3,\:n_3}}\left[\begin{array}{c}z_1\\z_2 \end{array} \left|\begin{array}{l} (a_j;\alpha_j,A_j;\xi_j)_{1,p_1} :(c_j,C_j;U_j)_{1,p_2}; (e_j,E_j;P_j)_{1,p_3}\\ (b_j; \beta_j, B_j; \eta_j)_{1,q_1} :(d_j,D_j ;V_j)_{2,q_2} ; (f_j,F_j;Q_j)_{1,q_3} \end{array} \right. \right] \eqno{(7.12)}$$
\normalsize
where $q_2 \geq n_2 \geq 1$, $ \Re(d)>0$.  \\\\
\textbf{(xiii)}\\
\small
$ \mathrm{{I}_{\:p_1,\:q_1 :\; p_2,\:q_2;\; p_3,\;q_3}^{\:0,\:n_1:\; m_2,\:n_2;\;m_3,\;n_3}} \left[\begin{array}{c}z_1\\z_2 \end{array} \left|\begin{array}{l} (a_j;\alpha_j,A_j;\xi_j)_{1,p_1} :(c_j,C_j;U_j)_{1,p_2}; (e_j,E_j;P_j)_{1,p_3}\\ (b_j; \beta_j, B_j; \eta_j)_{1,q_1} : (d_j,D_j ;V_j)_{1,q_2-1} ,(d,0;V); (f_j,F_j;Q_j)_{1,q_3} \end{array} \right. \right] $\\\\ 
$= \quad \frac{1}{\Gamma^{V}\left(1-d\right)} \quad \times $ \\
 $$\mathrm{{I}_{\:p_1,\:q_1 :\; p_2,\:q_2-1;\; p_3,\;q_3}^{\:0,\:n_1:\; m_2,\:n_2;\;m_3,\;n_3}} \left[\begin{array}{c}z_1\\z_2 \end{array} \left|\begin{array}{l} (a_j;\alpha_j,A_j;\xi_j)_{1,p_1} :(c_j,C_j;U_j)_{1,p_2}; (e_j,E_j;P_j)_{1,p_3}\\ (b_j; \beta_j, B_j; \eta_j)_{1,q_1} :(d_j,D_j ;V_j)_{1,q_2-1} ; (f_j,F_j;Q_j)_{1,q_3} \end{array} \right. \right] \eqno{(7.13)}$$
\normalsize
where $ q_2-1 \geq m_2 \geq 0$ , $ \Re(1-d) > 0 $ \\\\
\textbf{(xiv)}\\
\small
$   \mathrm{{I}_{\:p_1,\:q_1 :\; p_2,\:q_2;\; p_3,\:q_3}^{\:0,\:n_1:\; m_2,\:n_2;\;m_3,\:n_3}} \left[\begin{array}{c}z_1\\z_2 \end{array} \left|\begin{array}{l} (a_j;\alpha_j,A_j;\xi_j)_{1,p_1} :(c_j,C_j;U_j)_{1,p_2}; (e,0;P), (e_j,E_j;P_j)_{2,p_3}\\ (b_j; \beta_j, B_j; \eta_j)_{1,q_1} : (d_j,D_j ;V_j)_{1,q_2} ; (f_j,F_j;Q_j)_{1,q_3} \end{array} \right. \right] $\\ \\
$= \quad \Gamma^{P}\left(1-e\right) \quad \times $ \\
$$\mathrm{{I}_{\:p_1,\:q_1 :\; p_2,\:q_2;\; p_3-1,\:q_3}^{\:0,\:n_1:\; m_2\:n_2;\;m_3,\:n_3-1}} \left[\begin{array}{c}z_1\\z_2 \end{array} \left|\begin{array}{l} (a_j;\alpha_j,A_j;\xi_j)_{1,p_1} :(c_j,C_j;U_j)_{1,p_2}; (e_j,E_j;P_j)_{2,p_3}\\ (b_j; \beta_j, B_j; \eta_j)_{1,q_1} :(d_j,D_j ;V_j)_{1,q_2-1} ; (f_j,F_j;Q_j)_{1,q_3} \end{array} \right. \right] \eqno{(7.14)}$$
\normalsize
where $p_3 \geq n_3 \geq 1$, $ \Re(1-e)>0$. \\\\
\textbf{(xv)}\\
\small
$   \mathrm{{I}_{\:p_1,\:q_1 :\; p_2,\:q_2;\; p_3,\:q_3}^{\:0,\:n_1:\; m_2,\:n_2;\;m_3,\:n_3}} \left[\begin{array}{c}z_1\\z_2 \end{array} \left|\begin{array}{l} (a_j;\alpha_j,A_j;\xi_j)_{1,p_1} :(c_j,C_j;U_j)_{1,p_2}; (e_j,E_j;P_j)_{1,p_3-1},  (e,0;P)\\ (b_j; \beta_j, B_j; \eta_j)_{1,q_1} : (d_j,D_j ;V_j)_{1,q_2} ; (f_j,F_j;Q_j)_{1,q_3} \end{array} \right. \right] $\\\\ 
$=\quad \frac{1}{\Gamma^{P}\left(e\right)} \quad \times $\\
$$\mathrm{{I}_{\:p_1,\:q_1 :\; p_2,\:q_2;\; p_3-1,\:q_3}^{\:0,\:n_1:\; m_2,\:n_2;\;m_3,\:n_3}} \left[\begin{array}{c}z_1\\z_2 \end{array} \left|\begin{array}{l} (a_j;\alpha_j,A_j;\xi_j)_{1,p_1} :(c_j,C_j;U_j)_{1,p_2}; (e_j,E_j;P_j)_{1,p_3-1}\\ (b_j; \beta_j, B_j; \eta_j)_{1,q_1} :(d_j,D_j ;V_j)_{1,q_2-1} ; (f_j,F_j;Q_j)_{1,q_3} \end{array} \right. \right] \eqno{(7.15)}$$
\normalsize
where $p_3 - 1 \geq n_3 \geq 0$, $\Re(e)>0$.  \\\\
\textbf{(xvi)}\\
\small
$   \mathrm{{I}_{\:p_1,\:q_1 :\; p_2,\:q_2;\; p_3,\:q_3}^{\:0,\:n_1:\; m_2,\:n_2;\;m_3,\:n_3}} \left[\begin{array}{c}z_1\\z_2 \end{array} \left|\begin{array}{l} (a_j;\alpha_j,A_j;\xi_j)_{1,p_1} :(c_j,C_j;U_j)_{1,p_2}; (e_j,E_j;P_j)_{1,p_3}\\ (b_j; \beta_j, B_j; \eta_j)_{1,q_1} : (d_j,D_j ;V_j)_{1,q_2} ; (f,0;Q), (f_j,F_j;Q_j)_{2,q_3} \end{array} \right. \right] $\\ \\
$= \quad \Gamma^{Q}\left(f\right) \quad \times $ \\ $$ \mathrm{{I}_{\:p_1,\:q_1 :\; p_2,\:q_2;\; p_3,\:q_3-1}^{\:0,\:n_1:\; m_2,\:n_2;\;m_3-1,\:n_3}} \left[\begin{array}{c}z_1\\z_2 \end{array} \left|\begin{array}{l} (a_j;\alpha_j,A_j;\xi_j)_{1,p_1} :(c_j,C_j;U_j)_{1,p_2}; (e_j,E_j;P_j)_{1,p_3-1}\\ (b_j; \beta_j, B_j; \eta_j)_{1,q_1} :(d_j,D_j ;V_j)_{1,q_2} ; (f_j,F_j;Q_j)_{2,q_3} \end{array} \right. \right] \eqno{(7.16)}$$
\normalsize
where $q_3 \geq m_3 \geq 1$, $ \Re(f)>0 $. \\\\
\textbf{(xvii)}\\
\small
$   \mathrm{{I}_{\:p_1,\:q_1 :\; p_2,\:q_2;\; p_3,\:q_3}^{\:0,\:n_1:\; m_2,\:n_2;\;m_3,\:n_3}} \left[\begin{array}{c}z_1\\z_2 \end{array} \left|\begin{array}{l} (a_j;\alpha_j,A_j;\xi_j)_{1,p_1} :(c_j,C_j;U_j)_{1,p_2}; (e_j,E_j;P_j)_{1,p_3}\\ (b_j; \beta_j, B_j; \eta_j)_{1,q_1} : (d_j,D_j ;V_j)_{1,q_2} ;  (f_j,F_j;Q_j)_{1,q_3-1}, (f,0;Q) \end{array} \right. \right] $\\ \\
$ = \quad \frac{1}{\Gamma^{Q}\left(1-f\right)}\quad \times $\\
$$ \mathrm{{I}_{\:p_1,\:q_1 :\; p_2,\:q_2;\; p_3,\:q_3-1}^{\:0,\:n_1:\; m_2,\:n_2;\;m_3,\:n_3}} \left[\begin{array}{c}z_1\\z_2 \end{array} \left|\begin{array}{l} (a_j;\alpha_j,A_j;\xi_j)_{1,p_1} :(c_j,C_j;U_j)_{1,p_2}; (e_j,E_j;P_j)_{1,p_3-1}\\ (b_j; \beta_j, B_j; \eta_j)_{1,q_1} :(d_j,D_j ;V_j)_{1,q_2} ; (f_j,F_j;Q_j)_{1,q_3-1} \end{array} \right. \right] \eqno{(7.17)}$$
\normalsize
where $q_3-1 \geq n_3 \geq 0$ , $ \Re(1-f) > 0$. \\\\
\textbf{(xviii)}\\
\small 
$ \mathrm{{I}_{\:p_1,\:q_1 :\; p_2,\:q_2;\; p_3,\;q_3}^{\:0,\:n_1:\; m_2,\:n_2;\;m_3,\:n_3}} \left[\begin{array}{c}z_1\\z_2 \end{array} \left|\begin{array}{l} (a_j;\alpha_j,A_j;\xi_j)_{1,p_1} :\\ (b_j; \beta_j, B_j; \eta_j)_{1,q_1-1}, (a_1;\alpha_1,A_1;\xi_1) : \end{array} \right. \right. $ \\
$$ \qquad \qquad \qquad \qquad \qquad \qquad \left. \begin{array}{l} (c_j,C_j;U_j)_{1,p_2}; (e_j,E_j;P_j)_{1,p_3}\\ (d_j,D_j ;V_j)_{1,q_2-1} ,(c_1,C_1;U_1); (f_j,F_j;Q_j)_{1,q_3-1}, (e_1,E_1;P_1)\end{array}\right]$$
$ =  \mathrm{{I}_{\:p_1-1,\:q_1-1 :\; p_2-1,\:q_2-1;\; p_3-1,\:q_3-1}^{\:0,\:n_1-1:\; m_2,\:n_2-1;\;m_3,\:n_3-1}} \left[\begin{array}{c} z_1 \\ z_2 \end{array} \left|\begin{array}{l} (a_j;\alpha_j,A_j;\xi_j)_{2,p_1} : \\ (b_j; \beta_j, B_j; \eta_j)_{1,q_1-1}: \end{array} \right. \right. $  
$$ \qquad \qquad \qquad \qquad \qquad \qquad \left.\begin{array}{l}(c_j,C_j;U_j)_{2,p_2}; (e_j,E_j;P_j)_{2,p_3}\\ (d_j,D_j ;V_j)_{1,q_2-1} ; (f_j,F_j;Q_j)_{1,q_3-1} \end{array}\right] \eqno{(7.18)}$$ 
 \normalsize
provided that $p_1 \geq n_1 \geq 1$, $p_2 \geq n_2 \geq 1$, $p_3 \geq n_3 \geq 1$, $q_1 \geq 1$, $q_2 \geq m_2+1$ , $q_3 \geq m_3 +1$ \\ \\
\textbf{(xix)}\\
\small 
$ \mathrm{{I}_{\:p_1,\:q_1 :\; p_2,\:q_2;\; p_3,\:q_3}^{\:0,\:n_1: \;m_2,\:n_2;\;m_3,\:n_3}} \left[\begin{array}{c}z_1\\z_2 \end{array} \left|\begin{array}{l} (a_j;\alpha_j,A_j;\xi_j)_{1,p_1} :\\ (b_j; \beta_j, B_j; \eta_j)_{1,q_1} : \end{array} \right.\right. \qquad \qquad \qquad \qquad$\\
$$ \qquad \qquad \qquad \qquad \qquad \qquad  \left. \begin{array}{l}(c_j,C_j;U_j)_{1,p_2-1},(d_1,D_1;V_1); (e_j,E_j;P_j)_{1,p_3-1},(f_1,F_1;Q_1)\\(d_j,D_j ;V_j)_{1,q_2}; (f_j,F_j;Q_j)_{1,q_3} \end{array} \right] $$ 
$ =  \mathrm{{I}_{\:p_1,\:q_1 :\; p_2-1,\:q_2-1;\; p_3-1,\:q_3-1}^{\:0,\:n_1:\; m_2-1,\:n_2;\;m_3-1,\:n_3} }\left[\begin{array}{c}z_1\\z_2 \end{array} \left|\begin{array}{l} (a_j;\alpha_j,A_j;\xi_j)_{1,p_1} :\\ (b_j; \beta_j, B_j; \eta_j)_{1,q_1}: \end{array} \right. \right. $ 
$$\qquad \qquad \qquad \qquad \qquad \qquad \left. \begin{array}{l}(c_j,C_j;U_j)_{1,p_2 -1}; (e_j,E_j;P_j)_{1,p_3-1} \\(d_j,D_j ;V_j)_{2,q_2} ; (f_j,F_j;Q_j)_{2,q_3}  \end{array}\right] \eqno{(7.19)} $$
\normalsize
provided that $q_1 \geq 1$ , $q_2 \geq m_2 \geq 1$, $q_3 \geq m_3 \geq 1$ and $ p_1 \geq n_1+1$, $p_2 \geq n_2 +1$, $p_3 \geq n_3 + 1$.\\ \\
\textbf{(xx)}\\
\small 
$ \mathrm{{I}_{\:p_1,\:q_1 :\; p_2,\:q_2;\; p_3,\:q_3}^{\:0,\:n_1:\; m_2,\:n_2;\;m_3,\:n_3}} \left[\begin{array}{c}z_1\\z_2 \end{array} \left|\begin{array}{l} (a_j;\alpha_j,A_j;\xi_j)_{1,p_1} :\\ (b_j; \beta_j, B_j; \eta_j)_{1,q_1}  : \end{array} \right. \right. \qquad \qquad \qquad $ 
$$ \qquad \qquad \qquad \qquad \qquad \qquad \left. \begin{array}{l} (c_j,C_j;U_j)_{1,p_2}; \;(e_j,E_j;P_j)_{1,p_3}\\(c_{p_2},C_{p_2},U_{p_2}), (d_j,D_j ;V_j)_{2,q_2} ;\; (e_{p_3},E_{p_3};P_{p_3}), (f_j,F_j;Q_j)_{2,q_3}, \end{array}\right]$$
$$ =   \mathrm{{I}_{\:p_1,\:q_1 :\; p_2-1,\:q_2-1;\; p_3-1,\:q_3-1}^{\:0,\:n_1:\; m_2-1,\:n_2;\;m_3-1,\:n_3}} \left[\begin{array}{c}z_1\\z_2 \end{array} \left|\begin{array}{l} (a_j;\alpha_j,A_j;\xi_j)_{1,p_1} :\\ (b_j; \beta_j, B_j; \eta_j)_{1,q_1}  : \end{array} \right. \right. \qquad \qquad \qquad $$ 
$$ \qquad \qquad \qquad \qquad \qquad \qquad \left. \begin{array}{l} (c_j,C_j;U_j)_{1,p_2-1}; \;(e_j,E_j;P_j)_{1,p_3-1}\\ (d_j,D_j ;V_j)_{2,q_2} ;\; (f_j,F_j;Q_j)_{2,q_3}, \end{array}\right]  \eqno{(7.20)} $$ 
provided that   $ p_2 \geq n_2 + 1 $ , $p_3 \geq n_3 +1 $ , $q_2 \geq m_2$ and  $ q_3 \geq m_3$ \\\\
\textbf{(xxi)}\\
\small 
$ \mathrm{{I}_{\:p_1,\:q_1 :\; p_2,\:q_2;\; p_3,\:q_3}^{\:0,n_1-1:\; m_2,\:n_2;\;m_3,\:n_3}} \left[\begin{array}{c}z_1\\z_2 \end{array} \left|\begin{array}{l} (a_j;\alpha_j,A_j;\xi_j)_{2,p_1} :\\(b_{q_1},\beta_{q_1},B_{q_1}),  (b_j; \beta_j, B_j; \eta_j)_{1,q_1} : \end{array} \right.\right. $\\
$$ \qquad \qquad \qquad \qquad \qquad \qquad  \left. \begin{array}{l} (d_{q_2},D_{q_2},V_{q_2}),(c_j,C_j;U_j)_{2,p_2} ; (f_{q_3},F_{q_3},Q_{q_3}), (e_j,E_j;P_j)_{2,p_3}\\(d_j,D_j ;V_j)_{1,q_2}; (f_j,F_j;Q_j)_{1,q_3} \end{array} \right] $$
\\
$$  =     \mathrm{{I}_{\:p_1-1,\:q_1-1 :\; p_2-1,\:q_2-1;\; p_3-1,\:q_3-1}^{\:0,\:n_1:\; m_2,\:n_2-1;\;m_3,\:n_3-1}} \left[\begin{array}{c}z_1\\z_2 \end{array} \left|\begin{array}{l} (a_j;\alpha_j,A_j;\xi_j)_{2,p_1} :\\ (b_j; \beta_j, B_j; \eta_j)_{1,q_1-1}  : \end{array} \right. \right. \qquad \qquad \qquad $$ 
$$ \qquad \qquad \qquad \qquad \qquad \qquad \left. \begin{array}{l} (c_j,C_j;U_j)_{2,p_2}; \;(e_j,E_j;P_j)_{2,p_3}\\ (d_j,D_j ;V_j)_{1,q_2-1} ;\; (f_j,F_j;Q_j)_{1,q_3-1}, \end{array}\right]    \eqno{(7.21)} $$ 
where   $p_1-1 \geq n_1$ , $q_1 \geq 1$, $p_2 \geq n_2 \geq 1 $,  $p_3 \geq n_3 \geq 1 $ , $ q_2-1 \geq m_2$   and  $ q_3-1 \geq m_3$. \\\\ 
\normalsize
\textbf{Special Cases : } 
\par When all the exponents  $\xi_j (j=1,\ldots,p_1)$ , $\eta_j (j=1,\ldots,q_1)$ , $U_j (j=1,\ldots,p_2)$ ,  $V_j (j=1,\ldots,q_2)$ , $ P_j (j=1,\ldots,p_3)$ , $Q_j (j=1,\ldots,q_3)$  are  equal to unity , the I-function of two variables reduce to the H-function of two variables and therefore we obtain the \linebreak corresponding results in H- function of two variables recorded in [30].
\section*{Concluding Remark }
\hspace{20 pt}In this research paper we have introduced the natural generalization of the H-function of two variables introduced by Mittal and Gupta [21]. In addition to this we have also \linebreak obtained convergence conditions, series representation, elementary properties , transformation formulas and special cases.
\par As we have seen that the  generalized function of two variables introduced by \linebreak Agarwal [1] and Sharma [27,28] have found interesting applications in \textbf{wireless} \linebreak \textbf{communication}, therefore the function which have been  introduced  in this paper may be potentially useful.
\par We conclude  this paper by remarking that  single and double integrals, differentiation formulas and summation formulas  for this newly defined I-functions of two variables  are under investigations and the same  will form a subsequent paper in this direction.

  



\begin{thebibliography}{99}
\bibitem{1} Agarwal, R.P., \emph{An extension of Meijer's G-function}, Proc.Nat.Inst.Sci.India, Part A, {\bf 13}, 536-546, (1965)

\bibitem{2} Ansari, I.S., Yilmaz, F., Alouni, M.S. and Kucur, O.  \emph{ New results on the sum of Gamma random variates with application to the performance of Wireless communication \linebreak systems over Nakagami-m Fading Channels}, arxiv: 1202.2576v4 [cs.IT] 18 Jul. 2012
\bibitem{3} Ansari, I.S., Yilmaz, F.and  Alouni, M.S. ,  \emph{ On the sum of squared n-$\mu$ Random variates with application to the Performance of Wireless Communication Systems}, \linebreak  arxiv : 1210.0100v1 [cs.IT] 29, Sep 2012.
\bibitem{} Appel, P. and Kamp$\acute{e}$ de F$\acute{e}$riet, S., \emph{Fonctions hypergeometriques et hyperspheriques Polynomes d'Hermite}, Gauthier villars, Paris, (1926).   
\bibitem{4} Bailey W.N.,  \emph{ A Reducible Case of the Fourth Type of Appell's Hypergeometric Functions of Two Variables.} , Quart. J. Math. (Oxford) {\bf  4 }, 305-308, (1933).
\bibitem{7}Braaksma. B. L. J. , \emph{Asymptotic expansions and analytic continuations for a class of Barnes-integrals}, Compositio Math. \textbf{15}, 239-341, (1964).

\bibitem{5} Bora, S.L. and Kalla, S.L. , \emph{Some results involving generalized function of two variables} , Kyungpook Math. J. {\bf 10}, 133-140, (1970)
\bibitem{}  Burchnall. J.L. and Chaundy. T. W. \emph{Expansions of Appell's double hyper-geometric functions} , Quart. J. Math. , {\bf 11, no.1}, 249-270, (1940) 
\bibitem{}  Burchnall. J.L. and Chaundy. T. W. \emph{Expansions of Appell's double hyper-geometric functions (II)} , Quart. J. Math. , {\bf  12, no. 1}, 112-128, (1941)


\bibitem{8} Erdelyi.A. : \emph{Higher Transcendental Functions}, Vol.{\bf I}, McGraw-Hill Book Company, New York. (1953).
\bibitem{9} Fox,C. \emph{The G and H-functions as symmetrical Fourier kernals}, Trans.Amer. Math.Soc. {\bf 98}, 395-429, (1961).
\bibitem{10} Gupta.K.C., Jain Rashmi and Rajani Agrawal, \emph{On existence conditions for a \linebreak generalized Mellin - Barnes type integral}, Natl Acad Sci Lett, {\bf30(5,6)}, 169-172,  (2007).
\bibitem{11} Gupta.K.C. and Jain U.C., \emph{The H- function II}, Proc.Nat.Acad.Sci. {\bf36(A)}, 504-609, (1966).
\bibitem{12}Gupta.K.C.and Mittal P.K., \emph{Integrals involving a generalised function of two variables}, Indian J. Pure Appl. Math. {\bf 5}, 430-437, (1974).
\bibitem{} Goyal.S.P., \emph{The H-function of two variables} ,  Kyungpook Math.J. , {\bf 15, no.1}, 117-131,  (1975).
\bibitem{} Humbert, P. \emph{ The confluent hypergeometric functions of two variables}, Proc. Roy. Soc. Edinburgh Sect. A 41, 73-82, (1929).
\bibitem{13} Kilbas,A.A. and Saigo,M , \emph{H-transforms : Theory and applications} , Boca Raton - London - New  York- Washington , D.C. : CRC Press LLC, (2004). 
\bibitem{14}Luke.Y.L., \emph{The Special Functions and their Approximations}, Vol. 1, Academic Press, New York, (1969).

                                                                                                                                    
  



\bibitem{15}	Mathai.A.M., Saxena R.K., Haubold. H.J. :\emph{The H - function , Theory  and  Applications}, Springer, (2009).

\bibitem{17} Meijer, C.S. \emph { On the G-function I-VIII}, Nederl. Akad. Wetensch. Proc. {\bf 49}, 227-237; 344-356; 457-469; 632-641; 765-772; 936-943; 1063-1072; 1165-1175, (1946)
\bibitem{18} Mittal, P.K. and Gupta, K.C. : \emph{An integral involving generalized function of two \linebreak variables} , Proc.Indian Acad. Sci., sect.A. {\bf 75}, 117-123, (1972).



\bibitem{19} Pathak.R.S. \emph{Some results involving G and H functions} ,  Bull Calcutta. Math. Soc., 62, 97-106, (1970).
\bibitem{20} Rainville,E.D. \emph{Special Functions}, Macmillan Publishers, New York, (1963).
\bibitem{21}  Rathie, Arjun K.  \emph{A new generalization of generalized hypergeometric functions}, Le Matematiche  {\bf Vol. LII . Fasc. II},  297 - 310, (1997).
\bibitem{22} Ricard Beals, Roderick Wong   \emph{Special functions}, Cambridge University Press, (2010).




\bibitem{23} Sharma.B.L. and Abiodun.R.F.A. , \emph{Generating function for generalized function of two variables}, Proc.Amer.Math.Soc., {\bf 46 }, 69-72, (1974) 
\bibitem{}Sharma B L , \emph{ On the generalized function of  two variables (I)} , Annls.Soc.Sci.\linebreak Bruxeles Ser.{\bf I(79)}, 26-40 ,  (1965)

\bibitem{} Sharma , B.L. A new expansion formula for hypergeometric functions of two variables. Mathematical Proceedings of the Cambridge Philosophical Society, {\bf 64}, 413-416,  (1968) 
\bibitem{24} Srivastava.H.M. and Daoust. M. C., \emph{On Eulerian integrals associated with Kamp$\acute{e}$ de F$\acute{e}$riet function } , Publ. Inst. Math. (Beograd)(N. S.) , {\bf 9(23)} ,  199-202, (1969).
  
\bibitem{25} Srivastava.H.M., Gupta K.C., Goyal S.P., \emph{The H - functions of one and two variables with applications},South Asian Publishers, New Delhi.(1982).

\bibitem{26} Verma, R. U. , \emph{On the H-function of two variables I} ,  Indian J Pure Appl Math. {\bf5, No.7},  616-623 , (1974)

\bibitem{27} Xia Minghua, Wu Yik-Chung and Aissa Sonia, \emph{Exact Outage Probability of Dual-Hop CSI-Assisted AF Relaying Over Nakagami-m Fading Channels}, IEEE Transactions on Signal Processing, {\bf 60(10)}, 5578-5583, (2012)





\end{thebibliography}
\end{document}